\documentclass[12pt]{article}
\usepackage{amsmath, amsfonts, amssymb,amscd}

\def\l{\lambda}
\def\mfS{{\mathfrak S}}
\def\p{{\mathfrak p}}
\def\Sym{\mathfrak{Sym}}

\def\Part{\mathfrak{Part}}
\def\Car{\mathfrak{C\hspace{-0.07 em}a\hspace{-0.09 em}r}}
\def\Cos{\mathfrak{C\hspace{-0.07 em}o\hspace{-0.11 em}s}}
\def\Compo{\mathfrak{C\hspace{-0.07 em}o\hspace{-0.07 em}m\hspace{-0.09 em}p}}

\def\N{{\mathbb N}}
\def\Z{{\mathbb Z}}
\def\Q{{\mathbb Q}}

\def\C{{\mathbb C}}

\def\cH{{\mathcal H}}
\def\cC{{\mathcal C}}
\def\cN{{\mathcal N}}
\def\cS{{\mathcal S}}
\def\cY{{\mathcal Y}}
\def\hY{\widehat{\mathcal Y}}

\def\s{\scriptstyle}
\def\carre{\square\, }

\newtheorem{theorem}{Theorem}
\newtheorem{proposition}[theorem]{Proposition}
\newtheorem{lemma}[theorem]{Lemma}

\def\Proof{\noindent{\it Proof.}\ }

\def\CARRE#1{\hbox{\vrule width \thickness
   \vbox to \carresize{\hrule height \thickness\vss
      \hbox to \carresize{\hss#1\hss}
   \vss\hrule height\thickness}
\unskip\vrule width \thickness}
\kern-\thickness}
\def\vsquare#1{\vbox{\CARRE{$#1$}}\kern-\thickness}
\def\blk{\omit\hskip\carresize}

\def\smallyoung#1{%
  \newdimen\carresize \carresize=10pt%
  \newdimen\thickness \thickness=0.5pt%
  \vcenter{%
    \vbox{\smallskip\offinterlineskip%
      \halign{&\vsquare{##}\cr  #1}}}}

\def\young#1{%
  \newdimen\carresize \carresize=16pt%
  \newdimen\thickness \thickness=0.5pt%
  \vcenter{%
    \vbox{\smallskip\offinterlineskip%
      \halign{&\vsquare{##}\cr #1}}}}

\newdimen\unit
\def\o{$\scriptscriptstyle{{\rm o}}$}
\def\grape(#1,#2)#3{\raise#2\unit\rlap{\kern#1\unit #3}\ignorespaces}

\def\gg{{\unit=1mm
\hbox {\grape(2,1.2){'}
       \grape(1,2)\o
       \grape(2,2)\o
       \grape(3,2)\o
       \grape(1.5,1)\o
       \grape(2.5,1)\o
       \grape(2,0)\o\
}\kern 3.5 \unit}}

\def\gfill{\leaders\hbox to 1.2em{\hss\gg\hss}\hfill}

\def\frise{\medskip \centerline{\hbox to 8cm{\gfill} }\bigskip}
\begin{document}
\begin{center}
\bf\Large The Hecke algebra and structure constants\\
 of the ring of symmetric polynomials
\end{center}

\vspace{4mm}
\centerline{\large \emph{Alain Lascoux}}

\vspace{4mm}
\frise

\begin{abstract}
We give half a dozen bases of the Hecke algebra of the symmetric
group, and relate them to the basis of Geck-Rouquier, and to the basis
of Jones, using matrices of change of bases of the ring of
symmetric polynomials.
\end{abstract}


\bigskip
\section{Introduction}

To determine the characters of the symmetric group $\mfS_n$, 
Frobenius defined a linear morphism from $\C[\mfS_n]$, the group 
algebra of $\mfS_n$, to the space $\Sym$ of symmetric polynomials 
of degree $n$. This allowed him to identify the center of 
$\C[\mfS_n]$ and $\Sym$. Several natural bases of $\Sym$ 
correspond to specific bases of the center of $\C[\mfS_n]$.
This explains why transition matrices of symmetric functions 
occur as transition matrices of different bases of the center.

The identification between $\Sym$ and the center of the  Hecke algebra
of the symmetric group is not so clear. 
We still have a basis of central idempotents, and a basis,
due to Geck and Rouquier, which extends conjugacy classes.
 Jones \cite{Jones} (see also \cite{Francis-Jones})
 gave another natural basis. It happens that the transition
matrix of Jones' basis to the one of Geck and Rouquier
is, up to powers of a parameter $Q$, equal to the transition
matrix between power sums and monomial functions (Th. \ref{th:FJ}). 

This is this phenomenon that we want to explain in this text, 
and to generalize
by describing other bases in Theorems \ref{th:N1}, \ref{th:NTomega},
\ref{th:NCarre}. In our opinion, the most satisfactory 
explanation of the appearance of 
the transition matrices of the usual symmetric functions is due to 
the connection with
the theory of non-commutative symmetric functions 
given in Theorem \ref{th:NCSF}.

\bigskip
Recall that
the Hecke algebra $\cH_n$ of the symmetric group $\mfS_n$,
with coefficients in a commutative ring containing $q,q^{-1}$,
 is the algebra
generated by elements $T_1,\ldots, T_{n-1}$ satisfying the braid 
relations 
$$
\left\lbrace \begin{array}{c} T_iT_{i+1}T_i = T_{i+1}T_iT_{i+1} \\
    T_iT_j=T_jT_i \ ,\quad (\vert j-i \vert > 1) \end{array}\  , \right.  $$
together with the Hecke relations:
$$ (T_i-q)(T_i+1/q)=0  \, $$
or, equivalently, with $Q=q-1/q$,
\begin{equation}
 T_i^2= Q T_i +1  \, .
\end{equation}

Let us write $[k]$ for the $q$-integer  $(q^k-q^{-k})/(q-q^{-1})$,
and $[i.j\ldots k]$ for the product $[i]\, [j]\cdots [k]$.

To define the bases that we denote $\{ \cN_l(e_{\l^\natural}\}$,
$\{ \cN_\l(1) \}$, $\{ \cN_\l(T_{\omega_\l}^2 \}$,
 we need  the ring of coefficients  of $\cH_n$ to contain $Q^{-1}$
as well as the rational numbers. For the three other bases, we furthermore
need that the $q$-integers $[1],\ldots,[n]$ be invertible.

On $\cH_n$, one has a natural scalar product, with respect to
 which the basis
$\{T_w:\, w\in\mfS_n\}$ is orthonormal:
$$ (T_w, T_v)=  \delta_{w,v} \, . $$
It is such that 
$ (T_w T_i,T_v) = (T_w, T_v T_i)$, $ (T_iT_w,T_v) = (T_w, T_iT_v )$. 

Define $\Compo(n)$ to be the set of \emph{compositions} of $n$, i.e.
the set of vectors with positive integral components whose sum is $n$.
For any composition $I$, one defines a Young subgroup $\mfS_I$,
and its corresponding Hecke algebra $\cH(\mfS_I)$.
Elements of $\mfS_I$ are denoted $w^1\times w^2 \times \cdots$.
There is natural  
morphism, called
\emph{normalization}, from $ \cH(\mfS_I)$ to $\cH_n$:
\begin{equation}
\cH(\mfS_I) \ni h  \to  \cN_I(h):=  \sum_{w\in \mfS_n/\mfS_I}  
   T_w\, h\, T_{w^{-1}} \, .  
\end{equation}

One will find in \cite{Jones} many properties of more general norms.
We shall need only the 
following lemma, which shows that the normalization morphism 
allows to construct central elements.

\begin{lemma}
If $h$  is central in $\cH(\mfS_I)$, 
then  $\cN_I(h)$ is central in $\cH_n$.
\end{lemma}

Given $I$, one has also a projection  $\p^I: \cH_n \to \cH(\mfS_I)$,
which, for a permutation $w$, consists in considering it as a word cut
 into factors
of respective lengths $i_1,i_2,\ldots,i_r$, then renormalizing the 
values inside each factor.

We shall use different bases of symmetric functions,
and the corresponding matrices of change of basis $H2M, P2M, S2M,\ldots$
(see last section, and \cite{Mac}).

\section{Yang-Baxter basis}

Given an integral vector $v\in\Z^n$ (with components all different),
one defines a linear basis  $\{ \cY_w^v :\, w\in \mfS_n \}$ 
of $\cH_n$ recursively as follows:
$$ \cY_{ws_i}^v = \cY_{w}^v\, \left(T_i - \frac{q^k}{[k]}\right),\, 
        \ell(ws_i)>\ell(w),\, k=v_{w_{i+1}}-v_{w_i}\, ,$$
starting with $\cY_1=1$.

To any sequence of parameters $[z_1,\ldots,z_n]$ all different, 
one can in fact associate a Yang-Baxter basis. 
Here we have taken as ``spectral parameters''
 powers of $q$.  We shall use only the cases $v=[1,2,\ldots, n]$
and $v=[n,\ldots, 1]$. 

An important property of Yang-Baxter bases is a duality property
that is given in \cite{LLT} for general parameters. Taking into account
that we do not take the same scalar product, 
this duality formulates as follows.

% Let $\theta$ be the involution $q\to -1/q$.

\begin{theorem}   \label{th:yang}
% Let $\hY_w:= T_\omega\, \theta\left(\cY_{\omega w}  \right)$, $w\in\mfS_n$.
% Then $\{\hY_w\}$ is the basis adjoint to $\{\cY_w\}$, i.e.
%$$ \Bigl(\cY_v\, ,\, \hY_w  \Bigr)  = \delta_{v,w} \, .  $$
Given $v\in\Z^n$, let $u=[v_n,\ldots,v_1]$. Define 
 $\hY_w^v:= T_\omega\, \cY_{\omega w}^u$, $w\in\mfS_n$.
Then $\{\hY_w^v\}$ is the basis adjoint to $\{\cY_w^v\}$, i.e.
   $$ \Bigl(\cY_w^v\, ,\, \hY_{w'}^v  \Bigr)  = \delta_{w,w'} \, .  $$
\end{theorem}

The following two special cases are of interest, giving the two
1-dimensional idempotents of $\cH_n$, up to a normalizing factor 
(cf. \cite{LLT}).

\begin{lemma} 
  Let $\omega$ be the maximal permutation  
of $\mfS_n$.  Then 
\begin{eqnarray}   \label{nabla1}
\cY_{\omega}^{[1,\ldots,n]} 
        &=& \sum_{w\in \mfS_n} (-q)^{\ell(\omega)-\ell(w)}\, T_w \, ,\\
  \label{carre1}
\cY_{\omega}^{[n,\ldots,1]}          
      &=& \sum_{w\in \mfS_n} q^{\ell(w)-\ell(\omega)}\, T_w \, . 
\end{eqnarray}
\end{lemma}

\section{The center of the Hecke algebra}

It is clear that conjugacy classes
$\cC_\lambda:\, \lambda\in \Part(n)$ (considered as sums in the group
algebra) are a linear basis of the center of the group algebra of $\mfS_n$.

Geck and Rouquier \cite{Geck-Rouquier}
 have shown  that this basis extends canonically to
a basis $\Gamma_\lambda$ of the center of $\cH_n$.

According to Francis \cite{Francis}, 
the elements $\Gamma_\lambda$ are characterized by the 
property that
\begin{itemize}
\item $\Gamma_\lambda$ specializes to $\cC_\lambda$ for $Q=0$.
\item The difference $\Gamma_\lambda-\cC_\lambda$ involves 
no permutation which is of minimal length in its conjugacy class
(for the symmetric group).  
\end{itemize}

Let us write 
$$ \zeta_n = T_{n-1}\cdots T_2 T_1  \, ,$$
and define accordingly, by direct product,  for 
 any composition $I=[i_1,\ldots,i_r]$ of $n$, an element $\zeta_I$.
Explicitely,  let  
$K=[i_1,i_1+i_2,\ldots, i_1+i_2+\cdots +i_r]$. Then  
$\zeta_I$ is  the element $T_{w_I}$ indexed by
the following permutation in $\mfS_I$: 
$$ w_I= [k_1,1,\ldots,k_1-1,k_2,k+1+1,\ldots, k_2-1,\ldots, k_n,k_{n-1}+1,
   \ldots, k_m-1] \, .  $$

Reordering the factors, one sees that $\{  \zeta_I: I \in \Compo(n)\}$
is the set of subwords of $T_{n-1}\cdots T_2T_1$.

For example, for $I=[3,2,4]$, then 
$$ \zeta_I = \zeta_3 \times \zeta_2 \times \zeta_4 = 
  T_{312\, 54\, 9678}= (T_2T_1)\, (T_4)\, (T_8 T_7 T_6)= 
    (T_8 T_7T_6) (T_4) (T_2T_1)  \, .$$

There is one such $w_\lambda:\l\in\Part(n)$ 
in each conjugacy class for the symmetric group,
and it is of minimal length in its conjugacy class.
Therefore, any central element $g$  decomposes as the sum
$$ g =\sum_{\lambda\in\Part(n)}  
                 (g, \zeta_\lambda)\, \Gamma_\lambda \, .$$

We shall mostly use the following property of the elements $\zeta_J$, 
that is easy to check by induction.

\begin{lemma}  \label{th:ScalJM1}
  Given a composition $J$, given any permutation $w$,
then 
$$ (T_w,\, ,\,  \zeta_J T_w) = Q^{|J|-r} $$
if $w$ is of maximal length in its coset $w \mfS_J$,
and otherwise $ ( T_w\, ,\, \zeta_J T_w) =0$. 
\end{lemma}

This can be formulated in more striking terms. Define a \emph{recoil} of 
 a permutation $w$ to be any integer $i$ such that $i+1$ is left of $i$ 
 in $w$. Write $T_i\in\zeta_J$ if $T_i$ appears in a reduced decomposition
of $\zeta_J$. Then  $w$ is of maximal length in its coset $w \mfS_J$
iff $w$ has the recoil $i$ for all $i$ such that $T_i\in\zeta_J$.

An important remark follows from the lemma: the matrix 
representing the left multiplication by $\zeta_J$ in $\cH_n$ has the same
diagonal as the matrix representing the same multiplication 
in the $0$-Hecke algebra (with generators satisfying $T_i^2= QT_i$).
In other words, one can use the $0$-Hecke algebra to compute traces
of elements $ \zeta_J$.

\section{Jucys-Murphy elements}

Let the JM elements, first used by Bernstein,
 be $\xi_1=1$, $\xi_2= T_1T_1$, $\xi_3= T_2T_1T_1T_2$,
\begin{equation}
 \xi_i = T_{i-1} \cdots T_1\, T_1\cdots T_{i-1} \, . 
\end{equation}

These elements generate a subcommutative algebra (in fact maximal) of
$\cH_n$. The center of $\cH_n$ coincides with the symmetric functions
in $\xi_1,\ldots, \xi_n$ (this was a conjecture of 
Dipper-James \cite{DipperJames}, just settled by Francis and Graham
\cite{Francis-Graham}). 
However, these symmetric functions 
satisfy many relations which are not easy to control, 
and do not directly furnish the answer to our problem
of identifying the symmetric polynomials of degree $n$ with the
center of $\cH_n$. 

We shall write exponentially a monomial 
$$\xi_1^{u_1} \cdots \xi_n^{u_n} = \xi^u \, .$$
More generally, given any composition    $I=[i_1,\ldots, i_r]$, 
given any sequence of integral vectors $u\! u=[u^1,u^2,\ldots, u^r]$,
with $u^j\in \N^{i_j}$, let us write 
$$ \xi^{u\! u}  $$
for the direct product n $\cH(\mfS_I)$ of monomials in JM elements.
One can also view $u\! u$ as a single composition decomposed into factors of 
respective lengths $i_1,i_2,\ldots$ (we say that $u\! u$ is 
\emph{compatible} with $I$).

We shall only consider standard monomials $ \xi^{u\! u}$, i.e. such that
all components are boolean vectors (0-1 vectors). 

\begin{lemma}
Given a sequence $u\! u$ of boolean vectors, there exists a permutation
$\sigma(u\! u)$ such that 
$$ \xi^{u\! u} = T_{\sigma(u\! u)}\, T_{\sigma(u\! u)^{-1}} \, .$$
\end{lemma}

\Proof  The general statement is obtained by direct product from the
case of a single vector $u\in \{0,1\}^r$. 
Write
$$ \sqrt{\xi^u}  
:= (1)^{u_1} (T_1)^{u_2} \cdots (T_{r-1}\cdots T_1)^{u_r} \, . $$
This product is reduced, and therefore there exists a unique permutation 
$\sigma(u)$ such that $ \sqrt{\xi^u}=T_{\sigma(u)}$.
Moreover, using the commutation relations 
$ T_i\, (T_k\cdots T_1) = (T_k\cdots T_1)\, T_{i+1}$,
 $i<k$, one sees that 
$$ \xi^u =  \sqrt{\xi^u}  ( \sqrt{\xi^u})^\omega  \, ,$$
denoting by $z^\omega$ the reverse of
a word $z$.     \hfill QED

For example, for $u=[0,1,0,1,0]$,
$$ \sqrt{\xi^u} = T_1\cdot T_3T_2T_1 = T_{4213} \quad ,\quad
 \xi^u= T_1 T_3T_2T_1\cdot T_1T_2T_3T_1= T_{4213} T_{3241} . $$

Notice that for any composition $I$ and compatible $u\! u$, 
then all products $T_w T_{\sigma(u\! u)}$, for 
$w\in\mfS_n/\mfS_I$,  are reduced, and the permutations
$w \sigma(u\! u)$ are exactly those permutations  $w'$ such that
$\p^I(w') = \sigma(u\! u)$.

For example, with $I=[4,3]$ and $u\! u=[0101\, ,\, 001]$,
then the set 
$$\{ w \sigma(u\! u):\, w\in\mfS_n/\mfS_I\}  =
  \{ [4213\, 756],\ldots, [7546\, 312]\}   $$
consists of all permutations which project 
by $\p^{43}$ onto $[4213\, 756]$.  

\begin{theorem}   \label{th:ScalJM2}
Given two compositions $I,J$ of $n$, given $u\!u$
compatible with $I$, then
$$  ( \cN_I( \xi^{u\! u})\, ,\, \zeta_J)\,  Q^{\ell(J)-n}$$ 
is equal to the number of permutations in 
$\{ w\sigma(u\! u):\, w\in \mfS_n/\mfS_I \}$ having recoil $i$,   
for all $i$ such that $T_i\in \zeta_J$.
\end{theorem}

\Proof Extending the notation $\sqrt{ \xi^u}$ to the case of $u\! u$,
one rewrites
\begin{eqnarray*} 
 \left(\sum_{w\in\mfS_n/\mfS_I}T_w\, \xi^{u\! u}\,T_{w^{-1}}, \zeta_J\right)
 &=& \sum \left(T_w \sqrt{ \xi^{u\! u}}\, ,\,  \zeta_J\, T_w 
 \sqrt{ \xi^{u\! u}}\right)  \\
  &=&  \sum \left(T_{w\sigma(u\! u)}\, ,\,\zeta_J\,T_{w\sigma(u\! u)}\right)
\end{eqnarray*}
and one concludes with the help of Lemma \ref{th:ScalJM1}.   \hfill QED

For example, for $I=[3,2]$, $u\! u=[001,01]$, one has 
$\sqrt{\xi^{u\! u}}= T_2T_1\cdot T_4$, $\sigma(u\! u)=[3,1,2,5,4]$.
Taking $\zeta_J= T_1T_3$, one finds that there are two permutations 
in the set $\{ w\, [3,1,2,5,4]\}$ having recoils $1$ and $3$. 
They are $[42351]$ and $[52431]$. In consequence,
$$ (\cN_{32}( T_2T_1T_1T_2\, T_4T_4)\, ,\, T_1T_3) = 2\, Q^2 \, .$$

Taking the exponents $u\! u$ with all components equal to $0$, 
i.e. inducing from the identity, one gets a basis for the center:

\begin{theorem} \label{th:N1}
The set $\bigl\{ \cN_\l(1)\bigr\}_{\l\in \Part(n)}$  
is a basis of the center of $\cH_n$, and the matrix 
expressing it in the basis $\{ \Gamma_\l \}$ is
$$  E2M\cdot D\, ,$$
$D$ being the diagonal matrix with entries  
 $Q^{n-\ell(\l)}$, $\l\in \Part(n)$.
\end{theorem}

\Proof
According to the preceding theorem,
the entry $[\l,\mu]$  of the matrix is equal, up to powers of $Q$ to
the number of permutations $w\in \mfS_n/\mfS_\l \}$
having a recoil $i$ for all $i$ such that $T_i\in \zeta_\mu$.
But this the number of 0-1 matrices having row sums $\l$
and column sums $\mu$, which is also the coefficient of
the monomial function $m_\mu$ in the expansion of 
$e_\lambda= e_{\l_1} e_{\l_2}\cdots $.
Furthermore, the matrix having a non zero determinant,
the set $\{ \cN_\l(1)\}$ is a basis.    \hfill QED

For example, for $n=4$, the matrix is
$$  
\begin{bmatrix}
0   &   0   &   0   &   0   &   1   \\
0   &   0   &   0   &   Q   &   4   \\
0   &   0   &   Q^2   &   2Q   &   6   \\
0   &   Q^2   &   2Q^2   &   5Q   &   12   \\
Q^3   &   4Q^2   &   6Q^2   &   12Q   &   24   \\
\end{bmatrix}  \, .$$
Its third line is explained by the fact that, among the six 
permutations  in $\mfS_4/\mfS_{22}$, there are two with recoil $1$:
 $[2314]$ and $[2413]$, and one recoils $1,3$: 
$[2413]$.
% Induit1_2CC(4)   

\section{The Solomon module}

We shall see in this section that the square roots $\sqrt{\xi^{u\! u}}$
provide a connection with  non-commutative
symmetric functions.

A \emph{ribbon} $\theta$ is a skew diagram, not necessarily connected,
 (cf. \cite{Mac})
which contains no $2\times 2$ block of boxes.
One can write a ribbon as a sequence of compositions, recording the number
of boxes in each row, and passing to a new composition for every
connected component. For example,
$$  [ [3,1,2],[1,3]] = 
 \smallyoung{ & & \cr \blk &\blk & \cr \blk &\blk & & \cr 
 \blk &\blk &\blk &\blk & \cr \blk &\blk &\blk &\blk & & &\cr} \, .  $$

A permutation $w$ is \emph{compatible} with $\theta$ if, 
writing $w_1,w_2,\ldots$ in the successive boxes of $\theta$,
the result is a \emph{skew standard tableau}, 
i.e. increases in rows (from left to right), and columns
(from bootm to top).  A \emph{hook} is a connected ribbon of the type
$[1,\ldots,1,k]$.  

Given a ribbon $\theta$, the \emph{ribbon function} 
$R[\theta]$ is the sum, in the Hecke algebra, of all permutations
compatible with $\theta$. If $\theta$ is not connected, then 
$R[\theta]$ is equal\footnote{
This results from the iteration of the identity 
$$ R[ [\clubsuit,a],[b,\diamondsuit]] =   
   R[ [\clubsuit,a,b,\diamondsuit] + R[ [\clubsuit,a+b,\diamondsuit] \, , $$
where $\clubsuit,\diamondsuit$ 
are arbitray sequences, and $a,b$ are arbitrary positive integers,
cf. \cite{NCSF}. 
}
 to a sum of $R[J]$, $J$ composition of $n$.

The \emph{Solomon} module\footnote{As a subspace of $\C[\mfS_n]$,
it is a sub-algebra.} is the linear span of $R[J]$, $J\in \Compo(n)$.

We shall need another basis, generated by shuffle (for a general theory,
cf. \cite{NCSF} and \cite{JYT},
and also the different papers intitled {\tt NCSFx} on the page
of J-Y Thibon.  

Given a composition $J$, let $L[J]$ be the sum of all $T_w$: 
$w$ is  of maximal length in its coset $\mfS_J\, w$.
For example,
\begin{multline*}
 L[3,2] = T_{w_{3,2,1,5,4}}+T_{w_{3,2,5,1,4}}
            +T_{w_{3,2,5,4,1}}+T_{w_{3,5,2,1,4}}+T_{w_{ 3,5,2,4,1}} \\
+T_{w_{3,5,4,2,1}}+T_{w_{5,3,2,1,4}}+T_{w_{5,3,2,4,1}}
                    +T_{w_{5,3, 4,2,1}}+T_{w_{5,4,3,2,1}}\, . 
\end{multline*}

Notice that, for an y $w\in \mfS_n$, then
\begin{equation}  \label{ScalarL}   
 \bigl( T_w\, ,\, L[J]  \bigr) = 1 \ \text{or}\ 0   \ ,
\end{equation}
according to whether $w$ has a recoil in $i$ for all $i$ such that
$T_i\in \zeta_J$ or not.

Let $E_n(z)$ denote the product 
$$ (1+zT_1) (1+zT_2T_1)\cdots (1+zT_{n-1}\cdots T_1)  \, , $$
and, accordingly, define for a composition $J=[j_1,\ldots,j_r]$
the direct product 
$$ E_J(z_1,\ldots, z_r) 
           = E_{j_1}(z_1)\times \cdots \times E_{j_r}(z_r) \, .$$

Similarly, let $E_n^\xi(z)$ be the product
$$ E_n^\xi(z) := (1+z\xi_2) (1+z\xi_3)\cdots (1+z\xi_n)  \, , $$
and $E_J^\xi(z_1,\ldots, z_r)$ be the direct product 
$$ E_J^\xi(z_1,\ldots, z_r) = 
E_{j_1}^\xi(z_1)\times \cdots \times E_{j_r}^\xi(z_r) \, .$$

\begin{proposition}
For any composition $J$ of $n$, then 
$$ \sum_{w\in\mfS_n/\mfS_J} T_w\,  E_J(z_1,\ldots, z_r) $$
belongs to the Solomon module, being equal to
$$\sum_{I\leq J-1^r} z^I\, R\bigl[ [1^{i_1},j_1-i_1],\ldots,
  [1^{i_r},j_r-i_r] \bigr]\, .$$
\end{proposition}

\Proof  Every monomial appearing in the expansion of
$ E_J(z_1,\ldots, z_r)$ is reduced, and equal to some $\sigma(u\! u)$,
with $u\! u$ compatible with $J$. Conversely, one gets in this way all 
 $u\! u$ compatible with $J$
(the first component of each composition inside $u\! u$
has to be ignored, because $\xi_1=1$). 

In fact, $E_n(z)$ is equal to 
$\sum_{i=0^n-1} z^i R[1^i,n-i]$, and, by direct product,
$$ E_J(z_1,\ldots, z_r) = \sum_{I\leq J-1^r} z^I\, 
R\bigl[ [1^{i_1},j_1-i_1],\ldots, [1^{i_r},j_r-i_r] \bigr]\,
\bigcap \cH(\mfS_J) \, ,$$
each permutation belonging to $\mfS_J$ and being compatible
with a direct product of hooks. 

We have already noted that for any $w\in \mfS_n/\mfS_J$, the product
$T_w T_{\sigma(u\! u)}$ is reduced, and therefore 
the condition for the ribbons to belong to the subalgebra $\cH(\mfS_J)$
has to be lifted, multiplication by all $T_w$ permuting the values in such
a way as to preserve the ribbon shape. \hfill QED 

For example, for $J=[3,2]$, then 
\begin{multline*}
  E_{32}(z_1,z_2) = z^{00} R[[3],[2]]+ z^{10} R[[1,2],[2]]+
z^{01} R[[3],[11]]+ z^{20} R[[1,1,1],[2]] \\
+z^{11} R[[1,2],[1,1]] + z^{21} R[[1,1,1],[1,1]] 
\ \bigcap\, \cH_3 \times \cH_2 \, .
\end{multline*}
The terms 
$$ T_w\, \bigl((\sqrt{\xi_2}+\sqrt{\xi_3})\times \sqrt{\xi_2}  \bigr)
                 = T_w\,\bigl( (T_1+T_2T_1)\, T_4\bigr)  \, ,$$
$w\in \mfS_5/\mfS_{32}$, are exactly all the $T_v$ such that 
$v$ be compatible with $[[1,2],[1,1]]$, i.e. $v$ can be written as a 
ribbon tableau
$$ \young{v_1\cr v_2 &v_3\cr \blk &\blk &v_4\cr \blk &\blk &v_5 \cr} 
\ : \ v_1>v_2,\, v_2<v_3,\, v_4>v_5\, .  $$

\begin{theorem} \label{th:NCSF}
For any pair of compositions $J\in\N^r$, $K$ of $n$, then 
\begin{multline}
 \left(\sum\nolimits_{w\in\mfS_n/\mfS_J} T_w E_J(z_1,\ldots, z_r)\, ,\,  
 L[K]\right)   \\
 = \bigl(\cN_J(E^\xi_J(z_1,\ldots, z_r)\, , \, 
 \zeta_K \bigr)\, Q^{\ell(K)-n}  \, ,
\end{multline} 
and this polynomial in $z$ is equal to the coefficient of $m_K$ in the 
expansion of the product of symmetric functions
$$ \bigl(S_{1^{j_1}} +z_1 S_{2,1^{j_1-2}}+\cdots 
                       +z_1^{j_1-1} S_{j_1}  \bigr)\,\cdots   
\bigl(S_{1^{j_r}} +z_r S_{2,1^{j_r-2}}+\cdots +z_r^{j_r-1} S_{j_r} \bigr) 
\, .  $$
\end{theorem}

\Proof
Theorem \ref{th:ScalJM2} can be rewritten 
$$ \Bigl(\cN_J(\xi^{u\! u}\, ,\, \zeta_K  \Bigr) Q^{\ell(K)-n}  =
      \left(\sum\nolimits_{w\in\mfS_n/\mfS_J} 
           T_{w\sigma(u\! u)}\, ,\, L[K]\right) \, ,$$
because $L[K]$ is precisely the sum of all $T_v$ such that 
$v$ has a recoil in $i$ for all $i:\, T_i\in \zeta_k$.
Taking the generating function of all monomials $\xi^{u\! u}$, one gets 
$$ \Bigl(\cN_J(E_J^\xi(z_1,\ldots, z_r)))\, ,\, \zeta_K  \Bigr) Q^{\ell(K)-n}
  = \left(\sum\nolimits_{w\in\mfS_n/\mfS_J} T_w E_J(z_1,\ldots, z_r)
           \, ,\, L[K]\right) \, .$$
On the other hand, $\sum T_w E_J(z_1,\ldots, z_r)$ belongs to the 
Solomon module. Scalar products in this space can be evaluated 
at the commutative level (cf. \cite{Plactic}, \cite{NCSF}) and this 
gives the last statement.   \hfill QED

For example, for $J=[3,2]$, the scalar products 
\begin{multline*}
  \Bigl(\bigl( (1+z_1\xi_2)(1+ z_1\xi_3)\bigr) \times 
     \bigl(1+z_2\xi_2\bigr)\, ,\, L[K]\Bigr)   \\
 = \Bigl( (1+z_1T_1T_1)(1+z_1 T_2T_1T_1T_2)(1+z_2T_4T_4)\, ,\, L[K]\Bigr) 
\end{multline*}
and the scalar products 
$$  \Bigl(\cN_{32} \bigl( E_{32}^\xi(z_1,z_2)\bigr)\, ,\, \zeta_K\Bigr)
   Q^{\ell(K)-5} $$
are the coefficients of the monomial functions in the expansion of
\begin{multline*}
(S_{111} +z_1 S_{21} + z_1^2 S_3)\, (S_{11}+z_2 S_2)   \\
 = z^{21} m_5 + (z^{20}+z^{11}+2z^{21}) m_{41} +
  (z^{10}+2z^{11}+z^{20}+3z^{21}) m_{32}+\cdots  
\end{multline*}
To explain the coefficient, say, of $m_{41}$, we exhibit the contributing
permutations:
$$ z^{20}\, \smallyoung{\s 4\cr\s 3\cr \s 2\cr\blk &\s 1&\s 5\cr} 
+ z^{11}\, \smallyoung{\s 4\cr\s 3 &\s 5\cr \blk&\blk &\s 2\cr \blk&\blk &\s 1\cr} 
+ z^{21}\, \smallyoung{\s 4\cr\s 3\cr \s 2\cr\blk &\s 5\cr \blk &\s 1\cr}
 + z^{21}\, \smallyoung{\s 5\cr\s 4\cr\s  3\cr\blk &\s 2\cr \blk &\s 1\cr}\, . $$

\section{Basis of maximal products of JM elements}

Given a composition $I$, let us take 
$u\! u=u\! u(I):=[1^{i_1},\ldots, 1^{i_r}]$.
Then $\xi^{u\! u}$ is the product of all JM elements in each component
of $\cH(\mfS_I)$, and is equal to 
$T_{\omega_I}^2$, where $\omega_I$ is the permutation of maximal length
in $\mfS_I$.
 Moreover, the permutations $w \sigma(u\! u)$ 
are  maximal coset representatives for the cosets $\mfS_n/\mfS_I$
(i.e. are decreasing in each block).  

Taking the term of maximal degree in $z$ in Th. \ref{th:NCSF}, one gets:

\begin{theorem}   \label{th:NTomega}
The set of elements 
$\{\cN_\lambda\bigl(T_{\omega_\l}^2 \bigr),\, \lambda\in \Part(n)\}$,   
is a basis of the center of $\cH_n$.

The matrix expressing this basis in terms of the basis $\Gamma_\l$
is equal to
  $$   H2M\cdot D \, ,$$
$D$ being the diagonal matrix used in Th.\ref{th:N1}
\end{theorem}
% $$ D:= Diag\bigl( Q^{n-\ell(\lambda)} \bigr)_{\lambda\in \Part(n)}\, .  $$
% Base notee JJ 

For example, for $n=5$, 
writing the partitions in
the order $5, 41, 32, 311, 221$, $2111, 11111$  , the matrix 
of change of basis is equal to 
$$ \begin{bmatrix}
Q^4   &   Q^3   &   Q^3   &   Q^2   &   Q^2   &   Q   &   1   \\
Q^4   &   2 Q^3   &   2 Q^3   &   3 Q^2   &   3 Q^2   &   4 Q   &   5   \\
Q^4   &   2 Q^3   &   3 Q^3   &   4 Q^2   &   5 Q^2   &   7 Q   &   10   \\
Q^4   &   3 Q^3   &   4 Q^3   &   7 Q^2   &   8 Q^2   &   13 Q   &   20   \\
Q^4   &   3 Q^3   &   5 Q^3   &   8 Q^2   &   11 Q^2   &   18 Q   &   30
\\
Q^4   &   4 Q^3   &   7 Q^3   &   13 Q^2   &   18 Q^2   &   33 Q   &   60
\\
Q^4   &   5 Q^3   &   10 Q^3   &   20 Q^2   &   30 Q^2   &   60 Q   &   120
  \\
\end{bmatrix} \, .$$
The third line, for example, is obtained by expressing 
$\cN_{32}\left(T_1T_1 T_2T_1T_1T_2 T_4T_4\right) =
\cN_{32}\left(T_{32154}^2  \right)$ in the basis $\Gamma_\l$. 
% Projette(JJ([3,2]), 5);   

\section{Basis of Jones}

The first description of a basis of the center over $\Q[q,q^{-1}]$ 
was obtained by Jones \cite{Jones}. We shall rather follow the more 
recent paper of Francis and Jones \cite{Francis-Jones}, 
giving an explicit description
of this basis in terms of the basis of Geck and Rouquier.

We now use the original JM elements \cite{Jucys, Murphy}
 $x_1,\ldots,x_n$ defined 
by $x_1=0$, $x_2=T_1$, $x_3=T_2+T_2T_1T_2$,
 $$ x_j =\sum_{i<j} T_{(i,j)}   \, ,$$
sum over transpositions $(i,j)$.

It is immediate that 
$$ \xi_j= 1+Q x_j  \, ,$$
so that statements about the $\xi_j$'s can be translated in terms of
the $x_j$'s and conversely.

In $\cH_n$, let $e_i$, $i=0\ldots n-1$ denote the elementary symmetric
function of degree $i$ in $x_1,\ldots, x_n$. More generally,
given a composition $K=[k_1,k_2,\ldots]$,
let $e_{k_1}\times e_{k_2}\times \cdots$ 
denote a direct product of elementary
symmetric functions.

We have already induced central elements, starting from 
$1=e_0\times \cdots \times e_0$.
There is another case of special interest, taking the maximum 
(i.e non-vanishing) elementary symmetric functions. The following lemma
(due to Jones, \cite[Lemma 3.8]{Jones}) 
first describes the case of a composition with one part.

\begin{lemma}     \label{th:FJ2}
The maximal product $x_2\cdots x_n$ of JM elements is equal to 
$\Gamma_n$, and to 
$$ \cN_{1^{n-1}}(\zeta_n) := \sum_{w\in\mfS_{n-1}} T_w\, \zeta_n\,
  T_{w^{-1}} \, .$$
\end{lemma}

\Proof  Notice that 
$$Q^{n-1}\, x_2\cdots x_n=\xi_2\cdots \xi_n=(-1)^{n-1} E_{[n]}^\xi(-1)\, .$$
In the case $J=[n]$, Th. \ref{th:NCSF} states that, for any $K\in\Compo(n)$,
then
$$ \bigl((1-)^{n-1}E_{[n]}(-1)\, ,\, L[K] \bigr) =
\bigl((1-)^{n-1}E_{[n]}^\xi(-1)\, ,\, \zeta_K \bigr) $$
is the coefficient of $m_K$ in 
$$ (-1)^{n-1}\, \bigl(S_{1^n}- S_{2,1^{n-2}}+\cdots +
  (-1)^{n-1} S_n   \bigr) = m_n =p_n \, .$$
Therefore, only for $K=[n]$ does the scalar product 
$( x_2\cdots x_n\, ,\, \zeta_K)$ be different from $0$, and
$x_2\cdots x_n=\Gamma_n$.   

The element $X=\cN_{1^{n-1}}(\zeta_n)$ (considering 
$\mfS_{n-1}$ as the Young subgroup $\mfS_{n-1}\times \mfS_1$) 
is central in $\cH_n$. Each scalar product 
$$ \bigl( T_w\zeta_nT_{w^{-1}}\, ,\,  \zeta_J\bigr) =
   \bigl( T_w\zeta_n\, ,\, \zeta_JT_w \bigr)  $$
is null whenever $J\neq [n]$, because $T_w\zeta_n$ is equal
to some $T_v$, with $\ell(v)= \ell(w)+ n-1$, and no term of that 
length appears in the expansion of $\zeta_JT_w$. Since 
$$ T_w\, \zeta_n = \zeta_n T_{1\times w} \, ,$$
where $1\times w= [1,w_1+1,\ldots, w_{n-1}+1]$, the only permutation 
 $w$ which gives a non-zero scalar product is the identity.
\hfill QED   

By direct product, defining for any composition $I$ 
the composition  
$$ I^{\natural}:=[i_1-1,\ldots, i_r-1] \, $$
 then $e_{I^{\natural}}$ is the product of all non zero
JM elements in each component of $\cH(\mfS_I)$, as well as the direct
product $\Gamma_{i_1}\times \cdots \times\Gamma_{i_r}\times \cdots$. 

Specializing Th. \ref{th:NCSF} in $z_1=-1=z_2=\cdots$, and using   
Lemma \ref{th:FJ2}, one gets a basis of the center of $\cH_n$ 
\cite{Francis-Jones}
(Francis and Jones use another formulation, in term of the 
order\footnote{The matrix with entry $[\mu,\l]$,
$\mu,\l\in\Part(n)$, equal to 
 $$ n!\bigl|\cC_\mu \cap \mfS_\l \bigr|\, 
         \bigl|\cC_\mu\bigr|^{-1} \bigl|\mfS_\l \bigr|^{-1}$$ 
is in fact equal to $P2M$.
} 
of the intersections of conjugacy classes with Young subgroups).

\begin{theorem}   \label{th:FJ}
The set of elements 
$\bigl\{\cN_\lambda(e_{\l^\natural}) \bigr\}_{\lambda\in \Part(n)}$,
is a basis of the center of $\cH_n$.

The transition matrix to the basis $\{\Gamma_\l\}$
is equal to the product 
 $$ D^{-1}\cdot P2M \cdot D \, , $$ 
$D$ being the diagonal matrix used in Th.\ref{th:N1}.
%where $D$ is  the diagonal matrix
% with diagonal $[ Q^{n-\ell(\lambda)}:\lambda\in \Part(n)]$.   
\end{theorem}

For example, for $n=5$, writing the partitions in the order
$5, 41, 32$, $311, 221, 2111, 11111$, 
the matrix of Francis-Jones is 
$$ \begin{bmatrix}
1   &   0   &   0   &   0   &   0   &   0   &   0   \\
Q   &   1   &   0   &   0   &   0   &   0   &   0   \\
Q   &   0   &   1   &   0   &   0   &   0   &   0   \\
Q^2   &   2Q   &   Q   &   2   &   0   &   0   &   0   \\
Q^2   &   Q   &   2Q   &   0   &   2   &   0   &   0   \\
Q^3   &   3Q^2   &   4Q^2   &   6Q   &   6Q   &   6   &   0   \\
Q^4   &   5Q^3   &   10Q^3   &   20Q^2   &   30Q^2   &   60Q   &   120
  \\
\end{bmatrix} $$

One could have used other elementary symmetric functions in the 
JM elements. 
In the case where $I=[n]$, Lemma \ref{th:FJ2} generalizes to 
$$ e_k(x_2,\ldots, x_n) = \sum \Gamma_{J:\, \ell(J)=n-k}   \, ,$$
for any $k:\, 0\leq k <n$, and one can induce bases of the center
from direct products of such elements.

\section{Characters}

Central elements can also be determined by their characters.
In this section, we shall record properties of  characters
for the Hecke algebra, which will not be needed later,
but which illustrate the usefulness of both symmetric functions and 
Yang-Baxter elements and complete the relevant chapter of 
\cite{Geck-Pfeiffer}. 

In the case of the symmetric group $\mfS_n$, characters are constant 
on a conjugacy class, and consequently, it is sufficient to 
determine their values on a single element in each conjugacy class.
One way to compute characters is to use the 
\emph{Frobenius morphism} $\phi: \C[\mfS_n] \to \Sym$,
which sends a permutation $w$ to the product of power sums corresponding
to its cycle type decomposition:
$$ \phi(w)  = p^{cycle(w)}  \, .  $$
The characters of $w$ appear then as coefficients in the expansion 
of $\phi(w)$ in the basis of Schur functions. In fact, the 
Frobenius morphism is the easiest way to relate the center of
the algebra of the symmetric group to the ring of symmetric polynomials.

In the case of the Hecke algebra, one also call ``Table of characters''
the table of characters of the elements $\zeta_\lambda$ (which are canonical
elements of conjugacy classes), but a general element of $\mfS_n$ 
will not have its characters in this table.

The table of characters of $\cH_n$ has been determined by
Carter \cite{Carter}, Desarm\'enien \cite{Desarmenien}, Ram \cite{Ram}. 
Let us formulate their findings by transforming the Frobenius morphism.

For an alphabet $A$, write $\cS^k$ for the ``modified'' 
complete function
$$ \cS^k(A) := S^k(-Q A)\, Q^{-1},\ k \geq 1\, ,$$
and write exponentially $\cS^{ij\ldots}$ the product $\cS^i \cS^j\cdots$
of such functions.

Define $\Phi: \cH_n \to \Sym$ to be the morphism 
$$   \Phi(T_w) = \cS^{cycle(w)} \ .$$

Then, according to the cited authors, one has 

\begin{theorem}   \label{th:HekaCharTable}
Given a composition $J$ of $n$, the characters of $\zeta_J$ are the 
coefficients of the expansion of 
$(-1)^{\ell(J)} \cS^J$ in the basis of Schur functions.  
\end{theorem}

We shall now show that one can easily transform the table of characters
into a matrix independent of $q$.

Instead of $\zeta_J$, let us take, up to some normalizing
factor, the Yang-Baxter element corresponding to the same permutation.

Define 
$$ \Upsilon_n:= (T_{n-1}-q)\, \left(T_{n-2}-\frac{q^2}{[2]}  \right)
\cdots \left(T_1-\frac{q^{n-1}}{[n-1]}  \right)
 \left(\frac{-1}{[n]}  \right) \, $$
and, by direct product, $ \Upsilon_J$, for any composition $J$ of $n$.

For example, 
$$ \Upsilon_{43} = \left(T_3-q\right) \left(T_2-\frac{q^2}{[2]}  \right)
\left(T_1-\frac{q^3}{[3]}  \right) \left(\frac{-1}{[4]}  \right)
\left(T_6-q\right) \left(T_5-\frac{q^2}{[2]}  \right)
\left(\frac{-1}{[3]}  \right) \ .$$

The characters of such elements are simpler than those of $\zeta_J$.

\begin{theorem}  \label{th:CycleYang}
Given any composition $J$ of $n$, then 
$$ \Phi( \Upsilon_J)= h^J  \, .$$

The table of characters of all the $\Upsilon_\l: \l\in\Part(n)$
is the Kostka matrix $(E2S)^{tr}$.
\end{theorem}

\noindent{\it Proof.} The first statement comes, by direct product,
 from the case of a composition
with a single part. 
Let us show that one can compute $ \Phi( \Upsilon_n)$ by induction
on $n$,
filtering the product $(T_{n-1} -\alpha)\cdots (T_1-\beta)$ 
 according to the maximum right factor of type $T_j\cdots T_2T_1$.
For example, 
\begin{multline*}
 \left(T_3-q\right) \left(T_2-\frac{q^2}{[2]}  \right)
\left(T_1-\frac{q^3}{[3]}  \right) \left(\frac{-1}{[4]}  \right)  \\
=  T_3T_2T_1 \frac{-1}{[4]} + (-q) T_2T_1\frac{-1}{[4]} 
  + \left(T_3-q\right)\frac{-1}{[2]} q^2 T_1\frac{-1}{[4]}\\
 + \left(T_3-q\right) \left(T_2-\frac{q^2}{[2]}  \right)
 \frac{-1}{[3]} q^3 \frac{-1}{[4]} \, . 
\end{multline*}
Taking the image under $ \Phi$ of the right-hand side for a general $n$, 
one gets 
$\cS^n \frac{1}{[n]}+ \cS^{n-1}\, q S^1\frac{1}{[n]} +
\cdots + \cS^1\, q^{n-1}S^{n-1} \frac{1}{[n]}$. 

One notices then with pleasure that the first property 
enunciated in the theorem  reduces to
the identity $S^n( -QA+ qA)=S^n( -(q-1/q)A+ qA) = S^n(A/q)$
(see \cite{Cbms} for more identities using the $\l$-ring structure
of the ring of symmetric polynomials).      

The expansion of any $\Upsilon_J$ produces only subwords of
$T_{n-1}\cdots T_1$, which have the same characters as some
$\zeta_\lambda$. Therefore, the characters of $\Upsilon_J$ 
are the coefficients of the expansion of $\Phi(\Upsilon_J)$ in
the basis of Schur functions, and the matrix of characters of 
the $\Upsilon_\l$ is the transpose\footnote{We have here followed
the usual disposition for the table of characters, and this 
introduces a transposition.} of $E2S$ 
  \hfill  QED

Thus the table of characters\footnote{Using Th.\ref{th:CycleYang},
one can easily write determinants of elements of the Hecke algebra
which are sent under $\Phi$ onto Schur functions. The table of
characters now becomes the identity matrix.
For example, the determinant 
$\left|\begin{smallmatrix}   \Upsilon_3 & \Upsilon_4 \\ -1 &(T_4-q)(-1/[2])
\end{smallmatrix}\right|$, expanded by rows or by columns, 
is sent onto the Schur functions $s_{32}$. 
}
 is independent of $q$. Moreover, one could in fact
have directly computed the traces
of the matrices representing $\Upsilon_J$, and therefore obtained
the theory of characters of $\cH_n$ without using the results of
Carter, D\'esarm\'enien, Ram.

As an example, for $n=4$, the table of characters of $\zeta_4, \zeta_{31},
\zeta_{22}, \zeta_{211}, \zeta_{1111}$, written by columns, is 
$$   \Car_4 =  
\begin{bmatrix}
q^3   &   q^2   &   q^2   &   q   &   1   \\
-q   &   qQ   &   q^2-2   &   2q-1/q   &   3   \\
0   &   -1   &   [4]/[2]   &   Q   &   2   \\
1/q   &   -Q/q   &   1/q^2-2   &   q-2/q   &   3   \\
-1/q^3   &   1/q^2   &   1/q^2   &   -1/q   &   1   \\
\end{bmatrix}  $$
and the table of characters of $\Upsilon_4, \Upsilon_{31},
\Upsilon_{22}, \Upsilon_{211},\Upsilon_{1111}$ is
$$  \Car_4^\Upsilon =  \begin{bmatrix}
0   &   0   &   0   &   0   &   1   \\
0   &   0   &   0   &   1   &   3   \\
0   &   0   &   1   &   1   &   2   \\
0   &   1   &   1   &   2   &   3   \\
1   &   1   &   1   &   1   &   1   \\
\end{bmatrix} \ .$$

The first column of the last matrix tells us that
 $\Upsilon_4$ has zero trace in every irreducible
representation, except only one.

% The expansion of any $\Upsilon_J$ produces only subwords of 
%$T_{n-1}\cdots T_1$, which have the same characters as some 
% $\zeta_\lambda$.   

The matrix $\Car_n$ specializes, for $q=1$,
 to the table of characters $P2M^{tr}$ for the symmetric group 
(i.e. to the transpose of the matrix from the power sums
to the Schur basis).  

The relation between the two matrices has been studied by Carter 
\cite{Carter}, Ueno-Shibukawa \cite{Ueno}. 
The following property, that we shall not prove,
gives another factorization than theirs.

\begin{lemma}
For any $n$, let $D_3,D$ be the diagonal matrices with 
respective diagonals\footnote{
For a partition written exponentially: 
$\l=1^{\alpha_1} 2^{\alpha_2} 3^{\alpha_3}\cdots$, then 
$z_\l= 1^{\alpha_1} \alpha_1!\, 2^{\alpha_2} \alpha_2!
 \, 3^{\alpha_3} \alpha_3!\cdots $.
}
$ \bigl[ Q^{\ell(\lambda)-n}\prod_{i\in \lambda}[i] z_\l^{-1}
         ,\, \lambda\in \Part(n)\bigr]$,
$\bigl[ Q^{n-\ell(\lambda)},\, \lambda\in \Part(n)\bigr]$.
Then
$$ \Car_n = (P2S)^{tr}\cdot D_3\cdot P2M \cdot D \, .$$
\end{lemma}

For example, 
\begin{multline*}
 \Car_4= 
\begin{bmatrix}
1   &   1   &   1   &   1   &   1   \\
-1   &   0   &   -1   &   1   &   3   \\
0   &   -1   &   2   &   0   &   2   \\
1   &   0   &   -1   &   -1   &   3   \\
-1   &   1   &   1   &   -1   &   1   \\
\end{bmatrix}\cdot  
\begin{bmatrix}
\frac{[4]}{4}Q^{-3}   &   \cdot   &   \cdot   &   \cdot   &   \cdot   \\
\cdot   &  \frac{[3]}{3}Q^{-2}   &   \cdot   &   \cdot   &   \cdot   \\
\cdot   &   \cdot   &\frac{[2.2]}{8}Q^{-2}   &   \cdot   &   \cdot   \\
\cdot   &   \cdot   &   \cdot   & \frac{[2]}{4}Q^{-1}   &   \cdot   \\
\cdot   &   \cdot   &   \cdot   &   \cdot   &   \frac{1}{24}   \\
\end{bmatrix} \cdot  \\
\cdot 
 \begin{bmatrix}
1   &   0   &   0   &   0   &   0   \\
1   &   1   &   0   &   0   &   0   \\
1   &   0   &   2   &   0   &   0   \\
1   &   2   &   2   &   2   &   0   \\
1   &   4   &   6   &   12   &   24   \\
\end{bmatrix}   \cdot 
\begin{bmatrix}
Q^3   &   \cdot   &   \cdot   &   \cdot   &   \cdot   \\
\cdot   &   Q^2   &   \cdot   &   \cdot   &   \cdot   \\
\cdot   &   \cdot   &   Q^2   &   \cdot   &   \cdot   \\
\cdot   &   \cdot   &   \cdot   &   Q   &   \cdot   \\
\cdot   &   \cdot   &   \cdot   &   \cdot   &   1   \\
\end{bmatrix}
\end{multline*}

We have already met the matrix $P2M$ in the preceding section.

\section{Bases, by computing characters}

Let us now evaluate characters for different bases of central elements.

Recall that Young \cite{Young, Rutherford} 
defined orthogonal idempotents, 
indexed by standard tableaux,
for the group algebra of the symmetric group. These can be generalized 
to idempotents $\mho_t$ for $\cH_n$, the sum 
$$\mho_\l := \sum_t \mho_t $$ 
over all tableaux of shape $\l$ being a central idempotent.
Let us mention that one can view idempotents as polynomials
to compute them efficiently\cite{La97}.

The idempotents $\mho_t$ belong to the algebra generated by
the JM elements, and they are left and right eigenvectors
for them \cite{OV} :
$$ \mho_t \, \xi_i = \xi_i \, \mho_t = \mho_t\, q^{2c(i,t)} \, ,  $$
where $c(i,t)$ is the content of $i$ in the tableau $t$.

We shall not require these properties, but compute characters
by evaluating traces or scalar products.
Indeed, the value of the character of index $\l$ over an element $g\in\cH_n$
is equal to the scalar product    
$$ (\mho_\l\, ,\, g)   \, .$$

We first take the basis of  Jones, i.e. 
the family $\{ \cN_\lambda(  e^{\l^\natural}): \l\in\Part(n)  \}$.

\begin{theorem}   \label{th:FJChar}
Given $n$,
let $D_1$ and $D_2$  be the diagonal matrices with respective entries
 $$ \prod_i \frac{1}{[\l_i]}   \qquad ,\qquad 
 q^{\sum contents} n!\prod_{h=hook} \frac{[h]}{h}   \ , $$
product over all hooks $h$ of the diagram of $\l$,
for all $\l\in\Part(n)$.

Then the matrix of characters of the elements 
$\cN_\lambda\bigl(  e^{\l^\natural} \bigr)$ is equal to
$$ D_1\cdot P2S \cdot D_2  \ .    $$
\end{theorem}

For example, for $n=4$, the matrix of characters is
\begin{multline*}  
\begin{bmatrix}
q^6 [2.3]   &   -3q^2[2]   &   0   &   3q^{-2}[2] &   -q^{-6} [2.3] \\
q^6 [2.4]   &   0   &   -2[2.2]   &   0   &   q^{-6} [2.4]  \\
q^6 [3.4]/[2]   &-3q^2 [4]/[2] & 4[3] & -3q^{-2}[4]/[2] & q^{-6}[3.4]/[2]
  \\
q^6 [3.4]&   3q^2 [4]   &   0   &   -3 q^{-2} [4]   &   -q^{-6} [3.4] \\
q^6 [2.3.4]&   9q^2 [2.4] &   4 [2.2.3] &  9 [2.4] q^{-2} &   q^{-6}
[2.3.4]  \\
\end{bmatrix}  \\
 = \begin{bmatrix}
1/[4]   &   \cdot   &   \cdot   &   \cdot   &   \cdot   \\
\cdot   &   1/[3]   &   \cdot   &   \cdot   &   \cdot   \\
\cdot   &   \cdot   &   1/[2.2]   &   \cdot   &   \cdot   \\
\cdot   &   \cdot   &   \cdot   &   1/[2]   &   \cdot   \\
\cdot   &   \cdot   &   \cdot   &   \cdot   &   1   \\
\end{bmatrix}  
\cdot
\begin{bmatrix}
1   &   -1   &   0   &   1   &   -1   \\
1   &   0   &   -1   &   0   &   1   \\
1   &   -1   &   2   &   -1   &   1   \\
1   &   1   &   0   &   -1   &   -1   \\
1   &   3   &   2   &   3   &   1   \\
\end{bmatrix} \cdot \\  
\cdot \begin{bmatrix}
 q^6 [2.3.4]   &   \cdot   &   \cdot   &   \cdot   &   \cdot   \\
\cdot   &   3q^2 [2.4]   &   \cdot   &   \cdot   &   \cdot   \\
\cdot   &   \cdot   &   2 [2.2.3]   &   \cdot   &   \cdot   \\
\cdot   &   \cdot   &   \cdot   &   3 q^{-2} [2.4]   &   \cdot   \\
\cdot   &   \cdot   &   \cdot   &   \cdot   &   q^{-6} [2.3.4]   \\
\end{bmatrix}  
\end{multline*}

We are now in position to produce other bases of the center
by computing characters.

Given a composition $I$ define
$$ \carre_I :=\sum_{w\in \mfS_I} q^{\ell(w)} T_w   \quad \&\quad 
\nabla_I :=\sum_{w\in \mfS_I} (-q)^{-\ell(w)} T_w . $$

According to (\ref{nabla1}), (\ref{carre1}), 
both elements are proportional to some Yang-Baxter elements.
Indeed, 
$$ \carre_I= q^{\ell(\omega_I)}\, \cY_{\omega_I}^{[n,\ldots,1]}$$
and
$$ \nabla_I= (-q)^{-\ell(\omega_I)}\, \cY_{\omega_I}^{[1,\ldots,n]}\, ,$$
where $\omega_I$ is the permutation of maximal length in $\mfS_I$.

Both elements are central in $\cH(\mfS_I)$ and 
generate a 1-dimensional representation of $\cH(\mfS_I)$:
\begin{equation*}
  \carre_I T_w = \carre_I\,  q^{\ell(w)}
\quad \&\quad \nabla_I T_w= \nabla_I (-q)^{-\ell(w)}
\quad \text{for}\ w\in\mfS_I  \ .
\end{equation*}

% Let
% $$ \Pi_I = \cN_I\left( \carre_I  \right) =
%   \sum_{w\in\mfS_n/\mfS_I}  T_w\, \carre_I\,T_{w^{-1}} \, .$$

\begin{theorem} \label{th:NCarre}
Both sets $\bigl\{ \cN_\l(\carre_\l)\bigr\}_{\l\in\Part(n)}$ 
and $\bigl\{ \cN_\l(\nabla_\l)\bigr\}_{\l\in\Part(n)}$
are  a basis of the center
of $\cH_n$. 

Their characters  are respectively
$$   H2S\cdot D_2 \qquad \& \qquad E2S\cdot D_2 \ ,$$
where $D_2$ is the diagonal matrix defined in Th. \ref{th:FJChar}.
% $q^{\sum contents} n!\prod [h]/h$, product on all hooks of the diagram of $\l$.
\end{theorem}

\noindent{\it Proof.}  The fact that we have two bases will result from
the fact that their matrices of characters be invertible.

For each $\l$, the element 
$g_\l :=\sum_{w\in \mfS_n} T_w \carre_\l\,T_{w^{-1}} $ is equal to 
 $$\sum_{w\in\mfS_n/\mfS_\l,\, u\in \mfS_\l}  
   T_w\, T_u\, \carre_\l\, T_{u^{-1}}\, T_{w^{-1}}          
= \left(\sum_u q^{2\ell(u)}\right)\, \cN_\l(\carre_\l) \, .$$

Therefore, one can use $g_\l$ instead of $\cN_\l(\carre_\l)$.
However, the characters of  $g_\l$ are the same as the characters
of $\sum_{w\in \mfS_n} \carre_\l\,T_{w^{-1}} T_w$. 
Since $\sum_w T_{w^{-1}} T_w$ is a central element, 
we are finally reduce to compute the characters of 
$\carre_\l$.  

Given a representation of $\mfS_n$, let us decompose it into
irreducible subrepresentations of $\mfS_\l$. In the present case,
given a partition $\mu$, and the irreducible representation of
$\mfS_n$ with basis indexed by the standard tableaux of shape
$\mu$, we take any subspace with basis all the tableaux
differing from each other by a permutation in $\mfS_\l$.
Then $\carre_\l$ vanishes on this space, except if there is
a tableau $t$ containing the subwords $1,\ldots,\l_1$,
 $\l_1+1,\ldots,\l_1+\l_2$, $\ldots$,
i.e., which, considered as a permutation, is of minimal length in its
coset   $\mfS_\l\, t$. In that case, 
knowing that the trace of $\carre_k$ is equal to
$q^{k(k-1)/2}\, [1]\, [2]\cdots [k]$, then
the trace of 
 $\carre_\l$ is equal to 
$$ \prod_i  q^{i(i-1)/2} \prod_{j=1}^{\l_i} [j] \ . $$

The dependency on $\mu$ is therefore only the number of tableaux 
which are of  minimal length in their cosets  $\mfS_\l\, t$. 
This is one of the descriptions of the Kostka matrix $H2S$.

As for the set $\bigl\{ \cN_\l(\nabla_\l)\bigr\}$, one obtains it from
the first one by the involution $q\to -1/q$ and the exchange of the 
complete functions with the elementary ones (this corresponds to 
conjugating the partitions indexing irreducible representations).
\hfill QED

For example, for $n=4$, the matrix of characters of
$\cN_4(\carre_4), \cN_{31}(\carre_{31})$, $\cN_{22}(\carre_{22}), 
\cN_{211}(\carre_{211}), \cN_{1111}(\carre_{1111})$
(read by rows) is
$$ \begin{bmatrix}
1   &   0   &   0   &   0   &   0   \\
1   &   1   &   0   &   0   &   0   \\
1   &   1   &   1   &   0   &   0   \\
1   &   2   &   1   &   1   &   0   \\
1   &   3   &   2   &   3   &   1   \\
\end{bmatrix} \
\begin{bmatrix}
q^6[2.3.4]   &   \cdot   &   \cdot   &   \cdot   &   \cdot   \\
\cdot   &   3q^2[2.4]   &   \cdot   &   \cdot   &   \cdot   \\
\cdot   &   \cdot   &   2[2.2.3]   &   \cdot   &   \cdot   \\
\cdot   &   \cdot   &   \cdot   &   3q^{-2}[2.4]   &   \cdot   \\
\cdot   &   \cdot   &   \cdot   &   \cdot   &  q^{-6}[2.3.4]   \\
\end{bmatrix} \, ,  $$
the matrix of characters of
$\carre_4, \carre_{31},\carre_{22},\carre_{211}, \carre_{1111}$
is 
$$ 
\begin{bmatrix}
q^6[2.3.4]   &   \cdot   &   \cdot   &   \cdot   &   \cdot   \\
\cdot   &   q^3[2.3]   &   \cdot   &   \cdot   &   \cdot   \\
\cdot   &   \cdot   &   q^3[2.2]   &   \cdot   &   \cdot   \\
\cdot   &   \cdot   &   \cdot   &   q[2]   &   \cdot   \\
\cdot   &   \cdot   &   \cdot   &   \cdot   &   1   \\
\end{bmatrix}\ 
\begin{bmatrix}
1   &   0   &   0   &   0   &   0   \\
1   &   1   &   0   &   0   &   0   \\
1   &   1   &   1   &   0   &   0   \\
1   &   2   &   1   &   1   &   0   \\
1   &   3   &   2   &   3   &   1   \\
\end{bmatrix}     \, .
$$
the matrix of characters of
$\nabla_4, \nabla_{31},\nabla_{22},\nabla_{211}, \nabla_{1111}$
is
$$
\begin{bmatrix}
q^{-6}[2.3.4]   &   \cdot   &   \cdot   &   \cdot   &   \cdot   \\
\cdot   &   q^{-3}[2.3]   &   \cdot   &   \cdot   &   \cdot   \\
\cdot   &   \cdot   &   q^{-3}[2.2]   &   \cdot   &   \cdot   \\
\cdot   &   \cdot   &   \cdot   &   q^{-1}[2]   &   \cdot   \\
\cdot   &   \cdot   &   \cdot   &   \cdot   &   1   \\
\end{bmatrix}\
\begin{bmatrix}
0   &   0   &   0   &   0   &   1   \\
0   &   0   &   0   &   1   &   1   \\
0   &   0   &   1   &   1   &   1   \\
0   &   1   &   1   &   2   &   1   \\
1   &   3   &   2   &   3   &   1   \\
\end{bmatrix}
$$

The relations between the different bases that we have written 
may be summarized into the following diagram:

\begin{equation*}
\begin{CD} 
 \cN_\l(\nabla_\l) @>{E2P\cdot D_1^{-1}D^{-1}\cdot P2H}>> 
        \cN_\l(T_{\omega_\l}^2) @>{H2E}>> \cN_\l(1)\\ 
    @V{E2S\cdot D_2}VV  @VVV      @VV{E2P\cdot D}V    \\
 \mho_\l  @>>{D_2^{-1}\cdot S2H}>   
   \cN_\l(\carre_\l) @>>{H2P\cdot D_1^{-1}}> \cN_\l(e_{\l^\natural}) 
      @>>{D^{-1}\cdot P2M\cdot D}> \Gamma_\l 
\end{CD}
\end{equation*}

Ignoring the diagonal matrices, we have exhibited analogs 
of the bases $e,h,p,m,s$ of the ring of symmetric polynomials.
Notice that $\Gamma_\l$ should be identified with the basis of
monomial symmetric functions, and not with the basis of power sums.

In fact, we have produced two analogs of the bases $h^\l$ and $e^\l$. 
The matrix expressing the elements $\cN_\l(T_{\omega_\l}^2)$
in terms of the elements $\cN_\l(\carre_\l)$ is
$$ H2P\cdot D\cdot D_1\cdot P2H  $$
and becomes the identity when erasing $D$ and $D_1$
(but there is no specialization of $q$ which sends both $D$ and $D_1$
to the identity).

Moreover, it is remarkable that the matrix of change of basis 
between $\cN_\l\left(T_{\omega_\l}^2  \right)$   and 
$\cN_\l\left(1\right)$, which is $H2E$, be self-inverse and
independent of $q$.  Let us check  directly  this fact 
by showing how to relate 
normalizing $T_{\omega_\l}^2$ and normalizing $1$. 

Given a composition $J$ of $n$, consider the lattice of
compositions $I$ finer than $J$ (i.e. obtained by replacing the parts
of $J$ by compositions of them).

Let us write $\Cos(I)$ for the sum, in $\cH_n$, of the minimum coset
representative of $\mfS_n/\mfS_I$.  
It is easy to realize that
$$ \sum_{I\, finer\, J} (1-)^{n-\ell(I)} \Cos(I)  $$
is the sum of maximal coset
representative of $\mfS_n/\mfS_J$. 

Therefore
\begin{eqnarray*}
\sum_I \sum\nolimits_{w\in\mfS_n/\mfS_I}T_w\, T_{w^{-1}}  &=&
    \sum\nolimits_{w\in\mfS_n/\mfS_J}T_wT_{\omega_J} \, T_{\omega_J}  T_{w^{-1}} \\
   &=&  \cN_J\left( T_{\omega_J}^2  \right)  \, .
\end{eqnarray*}

On the other hand, the number of compositions $I$ finer than a
composition $J$, which reorder into a partition $\mu$ is,
indeed, the coefficient of $e^\l$ in the expansion of the product
$h^J$.  \hfill   QED

\section{Note: Symmetric Functions}
In this text, we have followed the conventions of Macdonald \cite{Mac},
and not of \cite{Cbms}, with the small proviso that we write
exponentially the multiplicative bases of the ring of symmetric 
polynomials: $e^\l$, $h^\l$, $p^\l$ are the products of elementary functions,
complete functions, power sums, respectively. 
One has also the monomial functions $m_\l$ and the Schur functions $S_\l$.

The matrix expressing the basis $a_\l$ (or $a^\l$) into the basis $b_\l$
(or $b^\l$) is denoted $A2B$, with capital letters. 
These matrices, for a given $n$,
 are obtained from the command \ {\tt Sf2TableMat(n,a,b)}\  in 
the maple library ACE \cite{ACE}.

The entries of such matrices may be written as scalar products,
and have many combinatorial interpretations. Realizing a Schur function,
or a skew Schur function as a sum of tableaux of a given (skew)-shape,
one gets in particular that the coefficient of $m_J$ in the
expansion of a product $h^I$ is equal to the number of 
non negative integral matrices
with row sums $I$, column sums $J$. Expanding  $e^I$ will 
similarly involve the 0-1 matrices with row sums $I$
and  column sums $J$ (see \cite[I. 6]{Mac} for more details).
Indeed, a product of monomials $x^u x^v\cdots x^w$ may be represented by
the matrix whose successive rows are $u,v,\ldots,w$. 

Consequently, a monomial $x^J$ appearing in a product $p^I$
corresponds to a matrix with row sums $I$ an column sums $J$, having 
only one non-zero entry in each row. This is another possible
description of the constants in
Th. \ref{th:FJ}.

\vspace{2cm}
\begin{center}
Alain Lascoux \\
CNRS, Institut Gaspard Monge, Universit\'e de Marne-la-Vall\'ee\\
77454 Marne-la-Vall\'ee Cedex, France\\[2pt]
{\tt Alain.Lascoux@univ-mlv.fr}
\end{center}
\end{document}